\documentclass[10pt]{amsart}
\usepackage{amssymb,amsmath,mathtools,enumitem}
\usepackage[T1]{fontenc}

\usepackage{mathabx}

\parskip=\smallskipamount

\newtheorem{theorem}{Theorem}%[section]

\newtheorem{corollary}[theorem]{Corollary}
\newtheorem{proposition}[theorem]{Proposition}

\theoremstyle{definition}
\newtheorem*{definition}{Definition}

\newtheorem*{examples}{Examples}

\newtheorem{remark}[theorem]{Remark}

\newcommand{\Z}{\mathbb{Z}}
\newcommand{\R}{\mathbb{R}}
\newcommand{\C}{\mathbb{C}}

\newcommand{\Gr}{\mathrm{Gr}}

\newcommand{\CP}{\mathbb{C}\mathrm{P}}

\newcommand{\G}{\mathrm{G}}
\renewcommand{\P}{\mathrm{P}}
\renewcommand{\c}{\mathrm{C}}

\numberwithin{equation}{section}

\begin{document}
%\doublespacing
\title[CR regular immersions and embeddings of $4$-manifolds into $\C^3$]
{CR regular embeddings and immersions of compact orientable $4$-manifolds into $\C^3$}
\author{Marko Slapar}
\address{University of Ljubljana, Faculty of Education, Kardeljeva Plo\v s\v cad 16, 1000 Ljubljana, Slovenia and Institute of Mathematics Physics and Mechanics, Jadranska 19, 1000 Ljubljana, Slovenia  }
\email{marko.slapar@pef.uni-lj.si}
\thanks{Supported by the research program P1-0291 and research project J1-5432 at the Slovenian Research Agency.}

%
%    General info
%
\subjclass[2000]{32V30, 32V40}
\date{\today} 
\keywords{CR manifolds, CR regular points, complex points, h-principle}

\begin{abstract}
We show that a compact orientable $4$-manifold $M$ has a CR regular immersion into $\C^3$ if and only if both its first Pontryagin class $p_1(M)$ and its Euler characteristic $\chi(M)$ vanish, and has a CR regular embedding into $\C^3$ if and only if in addition the second Stiefel-Whitney class $w_2(M)$ vanishes.  
\end{abstract}
\maketitle

\section{Introduction} 
All maps and manifolds in the paper are assumed to be smooth and by compact manifold we also mean closed (without boundary). Let $f\!:M\to X$ be an immersion of a real compact $2n$ dimensional manifold $M$ into a complex $n\!+\!1$ dimensional manifold $(X,J)$. A point $p\in M$ is called \textit{CR regular} for the immersion $f$ if the complex dimension of $f_*T_pM\cap Jf_*T_pM\subset T_{f(p)}X$ equals the expected dimension $n\!-\!1$. Points on $M$ that are not CR regular are called \textit{complex} or \textit{CR singular} points, and the complex tangent space at those points has the full dimension $n$. 
\begin{definition} An immersion (embedding) $f\!:M\to X$ is \textit{CR regular}, if all points on $M$ are CR regular for the immersion (embedding).
\end{definition}
Although it will also be obvious from the next section, we give a simple argument that gives a necessary condition for an oriented $4$-manifold $M$ to have a CR regular immersion into $\C^3$. Let $f\!:M\to \C^3$ be a CR regular immersion. Then $L=f^*(f_*TM\cap if_*TM)\to M$ is a complex line bundle on $M$. Let $E\to M$ be its orthogonal complement in $TM$ ($E$ is an $SO(2)=U(1)$ bundle). Since $(E\otimes\C)\oplus L=f^*T\C^3$ is a trivial complex bundle, it follows from the Whitney sum formula for Chern classes that
\[(1+c_1(E\otimes\C)+c_2(E\otimes\C))\cup (1+c_1(L))=1.\]
Since $c_1(E\otimes\C)=0$, we have $c_1(L)=0$, and thus also $c_2(E\otimes\C)=-p_1(E)=0$. So $L$ is trivial and, since $TM=L\oplus E$, $\chi(M)=0$ and $p_1(M)=p_1(E)=0$. We have shown that if $M$ can be CR regularly immersed in $\C^3$, is must satisfy $\chi(M)=0$ and $p_1(M)=0$. If, in addition, $f\!:M\to \C^3$ is an embedding, then by a result of Thom, the $\mathrm{mod}\ 2$ Euler class of the normal bundle of $M$ vanishes, and so by the Whitney sum formula, $w_2(M)$ also vanishes (we thank Peter Landweber for pointing this out). 

The main theorem of this paper shows that the above conditions are also sufficient for the existence of CR regular immersions and embeddings of compact oriented $4$-manifolds into $\C^3$. 

\begin{theorem}\label{thm} Let $M$ be a real compact orientable $4$-manifold. Then 
\begin{itemize}
\item[(i)] $M$ has a CR regular immersion into $\C^3$ if and only if $\chi(M)=0$ and $p_1(M)=0$,
\item[(ii)] $M$ has a CR regular embedding into $\C^3$ if and only if $\chi(M)=0$, $p_1(M)=0$ and $w_2(M)=0$. 
\end{itemize}
\end{theorem}
The vanishing of the second Stiefel-Whitney class $w_2(M)$ (and $M$ being orientable) is equivalent to $M$ being spin. By Hirzebruch's signature theorem we have $\sigma(M)=\frac{p_1(M)[M]}{3}$, where $\sigma(M)=b_2^+-b_2^-$ is the signature of (the intersection form of) $M$, so the vanishing of the first Pontryagin class $p_1(M)$ is equivalent to $\sigma(M)=0$. Since $SO(4)$-bundles are determined by their Euler, first Pontryagin and second Stiefel-Whitney class, we get the following corollary of the main theorem:

\begin{corollary} A compact orientable $4$-manifold has a CR regular embedding into $\C^3$ if and only if it is parallelizable. 
\end{corollary}

\begin{examples} We give some applications of the above theorem:
\begin{itemize}[leftmargin=*]
\item Since every compact orientable $3$-manifold $Y$ is parallelizable, $Y\times S^1$ can be CR regularly embedded into $\C^3$.
\item No simply connected compact orientable $4$-manifold can be CR regularly immersed into $\C^3$, since it always has the Euler characteristic at least $2$.
\item $\CP^2\#\overline\CP^2\#2(S^1\times S^3)$ can be regularly immersed into $\C^3$. In fact, any compact oriented $4$-manifold can be immersed into $\C^3$ after perhaps a connected sum with some copies of $\CP^2$ or $\overline\CP^2$, and $S^1\times S^3$ or $S^2\times S^2$.       
\end{itemize}
\end{examples}

\begin{remark} If $M$ is a compact oriented $4$-manifolds, then by a recent result of Jacobowitz and Landweber \cite{JL}, the vanishing of the first Pontryagin class $p_1(M)$ is equivalent to $M$ having both an independent map into $\C^3$ and a totally real immersion into $\C^4$.
\end{remark}  

\section{Complex points}

We assume from now on that $M$ is a compact oriented $4$-manifold, and $X$ is a complex $3$-manifold, and let $f\!:M\to X$ be an immersion. Let $Gf\!:M\to\widetilde{\Gr}_{4} f^*TX$ be the Gauss map of $M$ into the oriented real Grassmann bundle of the pull-back of $TX$, and let $(\P_\c f^*T^*X)^\pm\subset \widetilde{\Gr}_{4} f^*T^*X$ be the subbundle of complex hyperplanes, with their orientation either agreeing ($+$) or disagreeing ($-$) with the induced orientation from $TX$. By Thom's transversality theorem, for generic $f$ (which we also assume from now on), $\G f(M)$  and $(\P_\c f^*T^*X)^\pm$ intersect transversally, and so $M$ has only finitely many complex points by a simple dimension count. 

%\[\begin{tikzcd}[ampersand replacement=\&]
%(\P_\c f^*T^*X)^{\pm} \arrow[r,hook] \& \widetilde{\Gr}_{4} f^*TX\arrow[d] \arrow[r] \& \widetilde{\Gr}_{4} TX\arrow[d] \\
%\&M\arrow[u, bend left, "\widetilde{\G f}"] \arrow[r,"f"]\& X 
%\end{tikzcd}\]
A complex point $p\in M$ is called \textit{positive}, if at $p$, the manifold $M$ intersect  $(\P_\c f^*T^*X)^+$, and \textit{negative}, if it intersect $(\P_\c f^*T^*X)^-$; it is called \textit{elliptic}, if that intersection is positive, and \textit{hyperbolic}, if it is negative. Let $e_\pm(M)$ and $h_\pm(M)$ be the numbers of positive or negative elliptic or hyperbolic complex points on $M$. For \textit{Lai indices}, $I_\pm(M)=e_\pm(M)-h_\pm(M)$, we have topological formulas (\cite{Lai}) in terms of Chern classes of $X$ and Euler classes of the tangent and normal bundle of $M$:
\[2I_{\pm}(M)=(e(TM)\pm c_2(X)|_M+e(NM)\cup c_1(X)|_M\pm e(NM)^2)[M].\]
Since $e(NM)^2=p_1(NM)=p_1(TX)|_M-p_1(TM)$, in the case of $X=\C^3$, Lai formulas have a simple form
\[\label{lai2}2I_{\pm}(M)=\chi(M)\mp p_1(M).\]
 
Since Lai indices are invariant under regular homotopy of immersions, the vanishing of the indices is a necessary condition for the immersion to be regularly homotopic to a CR regular immersion. As a direct consequence of the Cancellation theorem (\cite{S1}) we have a converse statement: 

\begin{proposition}\label{cancel} Let $f\!:M\to X$ be a smooth generic immersion (embedding) of a real compact oriented $4$-manifold $M$ into a complex $3$-manifold $X$. Let $d$ be some metric on $X$ and let $\varepsilon>0$. If $I_\pm(M)=0$, there exists a regular homotopy (isotopy) $f_t:M\times [0,1]\to X$, so that                                                                       
\begin{itemize}
\item[(i)] $d(f_t(p),f(p))<\varepsilon$ for every $t\in [0,1]$ and every $p\in M$.
\item[(ii)] $f_1\!:M\to X$ is a CR regular immersion (embedding).
\end{itemize}                                     
\end{proposition} 

\begin{remark} The Cancellation theorem in \cite{S1} is proven for embeddings, but it is trivially adapted to immersions by passing to the normal bundle. The theorem holds in all dimensions and is essentially a consequence of Gromov's convex integration (\cite{Gr}). For real surfaces in complex surfaces, this is a classical result of Eliashberg and Harlamov (\cite{EH}).
\end{remark}

\section{Proof of Theorem \ref{thm}} 

The proof of the main theorem is a combination of classical results about immersions and embedding of orientable $4$-manifolds into $\R^6$ and Proposition \ref{cancel}. While some of the arguments below are probably well known to the experts, we include them both for completeness and because we were not able to find them in print. 

Let us first look at CR regular immersions into $\C^3$. We have already seen that the vanishing of the Euler characteristic and the first Pontryagin class is a necessary condition for the existence of CR regular immersions. Let us look at the converse statement. 

By Hirsh's immersion principle, an oriented compact $4$-manifold $M$ can be immersed into $\R^6$ if and only if there exists an $SO(2)$-bundle $L$ on $M$, so that $TM\oplus L$ is trivial ($L$ is the normal bundle of such an immersion). By a result of Cappell and Shaneson (\cite{CS}), the existence of such an $L$ is equivalent to the existence of an integral lift $\alpha\in H^2(M,\Z)$ of the second Stiefel-Whitney class $w_2(M)$, so that $\alpha\cup \alpha=-p_1(M)$. The cohomology conditions are obviously necessary. We sketch the proof of sufficiency (details can be found in \cite{K}). Let $L$ be an $SO(2)=U(1)$-bundle with $c_1(L)=\alpha$, $w_2(M)=\alpha\ (\mathrm{mod}\ 2)$ and $\alpha\cup \alpha=-p_1(M)$. Since $TX\oplus L$ is orientable, it can be trivialized over the $1$-skeleton of $M$. The obstruction to trivialization over the $2$-skeleton of $M$ is precisely $w_2(TX\oplus L)\in H^2(M,\pi_1(SO(6))=H^2(M,\Z_2)$, which by construction vanishes. Since $\pi_2(SO(6))$ is trivial, $TX\oplus L$ can be further trivialized over the $3$-skeleton of $M$, and the obstruction for the trivialization over the whole $M$ turns out to be $p_1(TX\oplus L)\in H^4(M,\pi_3(SO(6))=H^4(M,\Z)$. Since $p_1(TX\oplus L)=p_1(M)+c_1(L)^2=0$, $TX\oplus L$ is trivial.  

We now show that if $p_1(M)=0$, $M$ can be immersed into $\R^6$. By Hirzebruch's signature theorem, $\sigma(M)=0$, and by classification of signature $0$ unimodular forms (see for example \cite{GS}), the intersection form 
\[Q:H^2(M,\Z)/\mathrm{Tor}\times H^2(M,\Z)/\mathrm{Tor}\to\Z\] is equivalent to either $k\langle 1\rangle\oplus k\langle -1\rangle$ or $kH$, where $H=\left[\begin{smallmatrix} 0&1\\ 1&0\end{smallmatrix}\right]$. Since every oriented $4$-manifold is spin$^\mathrm{c}$, there exist an integral lift $\alpha\in H^2(M,\Z)$ of $w_2(M)$. Let $\alpha=\alpha_F+\alpha_T\in H^2(M,\Z)/\mathrm{Tor}\oplus \mathrm{Tor}$ be the splitting of $\alpha$ into the free and the torsion part.  Let us first assume that $Q$ is equivalent to $k\langle 1\rangle\oplus k\langle-1\rangle$, and let $a_1,a_2,\ldots,a_{2k}$ be the corresponding basis of $H^2(M,\Z)/\mathrm{Tor}$. Let us write $\alpha_F=n_1a_1+n_2a_2+\cdots+n_{2k}a_{2k}$ in this basis. By Wu's formula, $n_j=Q(\alpha_F,a_j)$ and $Q(a_j,a_j)$ have the same parity for every $1\le j\le 2k$, so $n_j$ must all be odd. Define $\widetilde\alpha=a_1+a_2+\cdots+a_{2k}+\alpha_T\in H^2(X,\Z).$ Since the $\mathrm{mod}\ 2$ reduction of $\widetilde \alpha-\alpha=(n_1-1)a_1+(n_2-1)a_2+\cdots+(n_{2k}-1)a_{2k}$ is trivial, $w_2(M)=\widetilde\alpha\ (\mathrm{mod}\ 2)$. We also have 
$(\widetilde\alpha\cup\widetilde\alpha)[M]=Q(a_1+a_2+\cdots+a_{2k},a_1+a_2+\cdots+a_{2k})=0=-p_1(M)[M],$ so $M$ can be immersed in $\R^6$. For $Q$ equivalent to $kH$, the proof is essentially the same, except that in this case, all $n_j$ must be even, and we can just take $\widetilde\alpha=\alpha_T$. 

Let now $M$ be an orientable manifold with $\chi(M)=0$ and $p_1(M)=0$. The above argument gives us an immersion of $M$ into $\C^3$ and since $2I_\pm(M)=\chi(M)\mp p_1(M)[M]=0$, we get that $M$ can be regularly homotoped to a CR regular immersion by Proposition \ref{cancel}.

We now prove the statement about embeddings. By a theorem of Cappell and Shaneson (\cite{CS}), an oriented $4$-manifold can be embedded into $\R^6$ if and only if $M$ is spin, $\chi(M)=0$ and $p_1(M)=0$ (for a complete proof, see \cite{R}). Again by Proposition \ref{cancel}, such an embedding can be isotoped to a CR regular embedding. 

\begin{remark} We can see from the proof above that if $\chi(M)=0$ and $p_1(M)=0$, one can construct many regularly non-homotopic CR immersions into $\C^3$ (unless $H^2(M,\Z)$ is completely torsion). If the intersection form of $M$ is equivalent to  $k\langle 1\rangle\oplus k\langle -1\rangle$, we could, in the above proof, take any $\tilde\alpha=m_1a_1+m_2a_2+\cdots+m_{2k}a_{2k}+\alpha_T$, as long as all $m_j$ are odd and $m_1^2+\cdots+m_k^2=m_{k+1}^2+\cdots+m_{2k}^2$. If the intersection form of $M$ is equivalent to  $kH$, we could take any $\tilde\alpha=m_1a_1+m_2a_2+\cdots+m_{2k}a_{2k}+\alpha_T$, as long as all $m_j$ are even and $m_1^2+m_3^2+\cdots+m_{2k-1}^2=m_2^2+m_4^2+\cdots+m_{2k}^2$. All these different choices provide us with regularly non-homotopic CR immersions, since they all have different normal bundles. 
\end{remark}

\begin{remark} For a final remark, notice that if a compact oriented $4$-manifold $M$ has a CR regular immersion into any complex $3$-manifold $(X,J)$, $f\!:M\to X$, then $TM$ splits as a direct sum of two $SO(2)=U(1)$-bundles, $TM=L\oplus E$, where $L=f^*f_*TM\cap f^*Jf_*TX$, and $E$ its orthogonal complement in $TM$ (in the case of $X=\C^3$, $L$ is trivial, as we have seen in the first section). So $M$ has an almost complex structure. A simple calculation using Lai's formulas shows, that all CR regular embeddings into $\CP^3$ are homologically trivial.    
\end{remark}       
      			                        
\bibliographystyle{amsplain}

\end{document}